\documentstyle[amscd,amssymb,verbatim]{amsart}

\theoremstyle{plain}
\newtheorem{Thm}{Theorem} 
\newtheorem{Lem}{Lemma}
\newtheorem{Prop}{Proposition}
\newtheorem{Def}{Definition}
\newtheorem{Con}{Conjecture}

\begin{document}

 
\title[Symplectic structure \\
on  the space of loops ]
{The structure of a symplectic \\
manifold on the space of \\
loops of 7-manifold }

\date{} 
\author{Michael Movshev} 
\address{Department
of Mathematics \\ 
SUNY at Stony Brook \\
 Stony Brook, NY 11794-3651 } 
\email{mmoushev@@math.sunysb.edu}

\dedicatory{}

\begin{abstract} We define a symplectic structure on the space of non parametrized loops in $G_2$ manifold. We also develop some basics of intersection theory of Lagrangian submanifolds. 
\end{abstract}

\maketitle

\section{Construction of symplectic structure}

\subsection{The transgression map}

There is a remarkable map described in \cite{Bry} which
relates de Rham complex on finite dimensional manifold $M$
with de Rham complex of its free loop space $\cal{L}(M)$:

\begin{equation}
 T:  \Omega^* ( M ) \rightarrow
\Omega^{*-1} ( \cal{L} (M) )
\end{equation}

The map $T$ is integration of a form along a path
$\gamma \in \cal{L}(M)$.  The main feature of this map
is that it is a complex morphism.  In particular it maps
closed forms into closed.  Observe that $Diff(S^1)^+$
acts on $\cal{L}(M)$ via reparametrization of a loop.
It is clear that the image of $T$ contains in
$Diff(S^1)^+$ invariant forms.  We want to employ $T$
to construct symplectic forms on the quotient $\cal{L}(M)$
by the action of the group.

\subsection{ Definition of the space of the singular knots}
Following \cite {Bry} we denote $\tilde X$ the subspace
in $\cal{L}(M)$ formed by immersions $\gamma$ which
have the properties that $\gamma$ induces an embedding
of $S^1 \backslash A$, for A finite subset of $S^1$ ,
and that for any two distinct points $x_1$ , $x_2$ of A
, the branches of $\gamma$ at $\gamma(x_1)$ and
$\gamma(x_2)$ have tangency of at most of finite order.
The group $Diff(S^1)^+$ acts freely on it.  We denote
the quotient by $\tilde Y$.

Following \cite {Bry} we identify a tangent space to a
singular knot $\gamma$ with the space of sections of
the normal bundle $N( \gamma)$ to $ \gamma$.

\subsection{The symplectic structure}

  Denote $ \lambda$ a closed 3-form on $M$.  $
\omega=T(\lambda)$ a closed 2-form on $ \tilde Y$.  To
be more explicit we give a formula for it.  Let $x,y$
be two tangent vectors to $ \gamma \in \tilde Y$.  We
identify them with sections of the normal bundle.  The
one form $ \lambda(.,x,y)$ is defined on $ \gamma$.
Then

\begin{equation} 
\omega (x,y) = \oint_{\gamma} \lambda
(*, x,y) 
\end{equation}

We state without a proof the following simple

\begin{Prop}
 The form $ \omega$ has no kernel iff $
\lambda$ satisfies some condition.  This condition is
purely local and must be fulfilled at any tangent
space.  Fix a tangent space $T_x$ to $x$.  The form $
\lambda(l, ., .)$ contains $l$ in its kernel for any
tangent vector $l$.The condition is that the kernel
must be one dimensional for any $l$.
\end{Prop}

The following lemma can be found in \cite {B}. We denote
by the same letter restriction of $\lambda$ to some
tangent space where it becomes a form with constant
coefficients.  The lemma gives a characterization of
possible $\lambda$.

\begin{Lem}
 1) The forms in question exist in
dimensions 3 and 7.

2) Let ${\Bbb A}$ be a noncommutative algebra with
division over ${\Bbb R}$.  The imaginary part
$Im({\Bbb A})$ of ${\Bbb A} $ is closed under commutator.
Let $(.  , .)$ be a natural dot-product on ${ \Bbb A}$.
Consider on $Im({\Bbb A})$ a form $([a,b],c)$, ($[a,b]$ is a commutator).  Then
$\lambda(a,b,c)$ is linearly equivalent to a multiple of
$([a,b],c)$.  
\end{Lem}

When dimension of the linear space is 3, $\lambda $ is
up to a constant a volume form.  In dimension 7 we
should start with Cayley numbers and follow the recipe
of the lemma.

It is clear that a space $ \tilde Y$ for 3-manifold,
equipped with a volume form bears a symplectic
structure.

The existence of form $\lambda$ on 7-manifold gives
restrictions on the structure group of the tangent
bundle.

\begin{Lem}
\cite {B}The group of automorphisms of
$\lambda$ in 7 dimensions is $G_2$.  
\end{Lem}

Now it should be clear how to construct such forms.
Pick a Riemannian 7 dimensional manifold with a holonomy
i $G \subseteq G_2$.  Then, making a parallel transport via
connection of a form defined over one tangent space, we
build a globally defined , closed form $\lambda$.

\subsection{Reconstruction of Cayley algebra from
tensor $\lambda$} \label{S:fgh}

According to \cite{B} all Cayley algebra can be
recovered from tensor $\lambda$.  Here a short sketch
how he does it.

We know that Cayley numbers ${\Bbb O}$ are equipped
with canonical anti-involution , trace and a dot
product.  By definition $Im {\Bbb O}=\left \{ l \in
{\Bbb O}; \bar{l} =-l \right \}$.  Here $\bar .$ is
antiinvolution.  Then $tr(l)=1/2(l+\bar l)$,  $
(l,m)=tr(l* \bar m)$ and $\lambda(a,b,c)=([a,b],c)$.
To recover $[.,.]$ from $\lambda$ we have to know  a dot product.  Define a
linear map ${\Bbb R}^7 \otimes {\Bbb R}^7 \rightarrow
\Lambda^7(({\Bbb R}^7)^* )$by the formula $m \otimes n
\rightarrow \lambda (n,.,.) \wedge \lambda (n,.,.)
\wedge \lambda$.  This defines a conformal class of
a metric.  To fix a concrete representative we chose
it satisfying 

\begin{equation} \label{E:gty}
[m,[m,n]]=-(m,m)n+(m,n)m 
\end{equation}

\subsection{A simple observation} \label{S:asw}

      Denote as usual $TM$ a tangent bundle of a
manifold $M$.  Denote the compliment of the zero
section in $TM$ by $T^0(M)$.  Let $p: T^0(M)
\rightarrow M$ be a natural projection.  Let
$p^*T$ be the pullback of the tangent bundle.  It
contains a tautological subbundle $O$.  We denote
by $I$ the quotient $p^*T/O$.  Every element
$\gamma$ of $\tilde X$ defines a map $
\overset{.}{\gamma}:  S^1 \rightarrow T^0(M)$.
We can identify canonically the tangent space to
$\gamma$ in $\tilde Y$ with sections of the
pullback $\gamma^*I$.

It is clear that $I$ in case of 3 dimensional
Riemannian manifold and $G_2$ manifold is a
symplectic bundle.  The symplectic stricture
inherits from the 3-form.  More precisely we
pullback the form $\lambda$ on $T^0(M)$.  $T^0(M)$
admits a nonvanishing vector field $e$:  a
generator of dilatations along the fibers of the
projection.  This field is a basis of $O$.  We
plug $e$ into $\lambda$.  The result is a
nondegenerate form $\omega$ on $I$.

Acting the same way we can define a complex
structure on $I$.  In 3 dimensional case we pick
a Riemannian metric, whose volume form is
$\lambda$.  In 7 dimensions according to
\ref{S:fgh} we have some canonical Riemannian
metric.  In both cases we have a vector
product $[.,.]$ at hand.  Denote the normalization
of $e$ by  $e_0$ and identify $I$ with orthogonal
compliment to $e$.  Define a complex structure on
$I$ via the formula $J(.)=[e_0,.]$. The identity $J^2=-Id$ is a
corollary of \ref{E:gty} and a similar identity, 
which holds in 3  dimensions.

We want to bring to your attention a fact that symplectic,
Riemannian and complex structure on $I$ are compatible in
 a sense that $\omega(.,.)=(J.,.)$

\section{Some topological computations}
Let $I$ be a complex bundle, defined in \ref{S:asw}.

\begin{Lem} \label{L:fvx} 
The first Chern class of $I$, defined over the base $T^0(M)$
($M$ is 3 dimensional) is divisible by 2.
\end{Lem} 
 
\begin{pf}
 We use the fact that the tangent
bundle of an oriented 3 manifold is parallelizable. We replace 
$T^0M$ by the a homotopy equivalent to it spherical bundle $S(M)$.
 The spherical bundle $S(M)$ is trivial. Denote trivializing
projection $S(M) \rightarrow S^2$ by $q$ The bundle $I$ is
isomorphic to the pullback  $q^*(T(S^2))$ of the tangent bundle $T(S^2)$ of  $S^2$. This being said the claim about Chern class becomes trivial.
 \end{pf}
  
We want to have a similar information about $G_2$
manifold $M^7$.  The spherical bundle over $M^7$ also
has a complex bundle $I$ on it. $I$ again can be 
identified  with a subbundle  of tangent bundle of 
$S(M)$ consisting of vectors tangent to fibers
 of projection $p:S(M) \rightarrow M$. 

 \begin{Lem} \label{L:wsd}
The first Chern class of $I$ over $T^0(M^7)$ is zero.
  \end{Lem}
  
\begin{pf}

As before we replace $T^0(M^7)$ by the spherical 
bundle $S(M)$. We already mentioned that  $I$ is the tangent bundle to
the fibers of projection $p:S(M) \rightarrow M$.  We
want to show  that the top wedge power $\Lambda_{{\Bbb C}}^3(I)$ has a nonvanishing 
section.  

A sphere  $S^6$ is  $G_2$ homogeneous space. with stabilizer
 of a point being equal to $SU(3)$.
$\Lambda_{{\Bbb C}}^3(T(S^6))$ must contain a  $G_2$ 
invariant nonzero section. It gives rise to a never vanishing 
section of $\Lambda_{{\Bbb C}}^3(I)$.
\end{pf}

\section{Lagrangian submanifolds in $\tilde Y$}

In \cite {Bry} some examples of Lagrangian
submanifolds were proposed.  Namely take $L{\Sigma}$
be a set of all loops containing in the embedded two
dimensional surface $\Sigma \subset M^3$. In \cite{Bry}
it was proved that it is a Lagrangian submanifold of
 $\tilde Y$.

In $G_2$ setting the role of hypersurfaces is played by
isotropic submanifolds, in other words embedded 4-manifolds ,
the restriction of  form $\lambda$ on whose is  equal to
zero identically.

 \subsection{Isotropic subspaces} \label{S:ert}

We will need some information on linear subspaces of
${\Bbb R}^7$, the restriction of tensor $\lambda$ on
whose vanishes identically.  We call such subspace
isotropic.

\begin{Lem} 

1) The dimension of maximal isotropic
subspace is 4.

2) The group $G_2$ acts transitively on the set of
isotropic subspaces.

3) The stabilizer of isotropic 4-plane is $SO(4)$.

\end{Lem}

\begin{pf}

Without loss of generality we assume that
$\lambda(a,b,c)=([a,b],c)$, where commutator and dot
product are taken  in Cayley numbers.  Take a
nonzero vector from a maximal isotropic subspace $L$
and call it $\bold i$.  Take any vector in orthogonal
compliment to $\bold i$ and call it $\bold j$.
Obviously $[\bold i,\bold j] \perp L$, otherwise the
restriction of $\lambda $ on $L$ would be nonzero.  Let
us call $k$ the Cayley vector product $[\bold i, \bold
j]$.  Pick any vector in the orthogonal compliment to
$\bold i, \bold j $ in $L$.  We call it $\bold l$.  It
is clear that $\lambda|_{< \bold i, \bold j, \bold l >
} \equiv 0$.  We are claiming that $[\bold k, \bold
l] \in L$ because $L$ is maximal isotropic.  Obviously
$\lambda$ remains to be identically zero on the space
$<\bold i, \bold j, \bold l,[\bold k, \bold l ]>$.  If
$L$ contained some other vector than $[\bold k, \bold
l]$, orthogonal to $<\bold i, \bold j, \bold l>$, then
by simple reasonings $\lambda|_L \neq 0$.  It is also
clear that we can add no more vectors to enlarge $L$.  This
proves first assertion.

If we normalize vectors $\{\bold 1, \bold i, \bold j,
\bold k,\bold l, [\bold i, \bold l],[\bold j, \bold
l],[\bold k, \bold l]\}$, then in each such a basis
multiplication in ${\Bbb O}$ is given by the same
structure constants.  A linear transformation that
carries one such a basis into another is an automorphism
of ${\Bbb O}$ and hence is in $G_2$.  The above basis is
a function of $\bold i, \bold j, \bold l$.  This implies
that $G_2$ acts transitively on  triples $\bold i, \bold j, \bold l$ and therefor on isotropic subspaces.  This proves the second
assertion.

   It is clear that the stabilizer St(L) acts transitively and free on
the space of triples of orthogonal vectors in $L$.  It
identifies $St(L)$ with $SO(4)$.
\end{pf}

\begin{Prop} 
Let $\Sigma$ be immersed isotropic
submanifold.  The set $L(\Sigma)$ of loops
containing in $\Sigma$ is a Lagrangian submanifold.
\end{Prop}

\begin{pf}

 It is clear that the restriction of
symplectic form $\omega$ on $L(\Sigma)$ is zero.  Choose a
parametrization of a curve $\gamma \in \tilde Y$.  The
tangent space to $\gamma$ is identified with sections
of normal bundle to $\gamma$.  The subspace of tangent
space , tangent to $L(\Sigma)$ is the space of normal
bundle in $\Sigma$.  Then it is clear that any vector ,
orthogonal in a sense of $\omega$ to sections of normal
bundle in $L(\Sigma)$ must be contained in it. 
\end{pf}
Now we want to study how many connected components $L(\Sigma)$ has.

Following \cite{Bry} we define a space $Map(S^1,S(\Sigma))$. $S(\Sigma)$ 
is a spherical bundle associated with the tangent bundle of 
$\Sigma$. $S^1$ acts on $Map(S^1,S(\Sigma))$. 
\begin{Prop}
$L(\Sigma)$ is weakly homotopy equivalent to the Borel construction
 of $Map(S^1,S(\Sigma))$.
\end{Prop}

This enables us to compute the number of connected components of $L(\Sigma)$.

 We stick for a moment to the case of  surfaces in 3 dimensional manifold.
We state without the proof 
\begin{Prop}
The connected components of $L(\Sigma)$ are in one to one correspondence
 with elements of $\pi_1(S(\Sigma))$. It is ${\Bbb Z}_2$ in case $\Sigma=S^2$ and 
 canonical central extension 
 $1 \rightarrow {\Bbb Z}\rightarrow \pi_1(S(\Sigma)) \rightarrow \pi_1(\Sigma) \rightarrow 1$ of the
  fundamental group of a surface $\Sigma$.
\end{Prop}

Observe that since we do care about orientation of curves in $L(\Sigma)$,
 we have two curves corresponding the underlying set of points.

\subsection{Some Lagrangian subbundles of  $I$} \label{S:tgu}
In \ref{S:asw} we defined a symplectic, complex vector bundle over $T^0(M)$
 the compliment to the zero section of the tangent bundle.
 
 Now we suppose that we have an  
 embedding  $i:\Sigma\rightarrow M$,
   $\Sigma$ being a two dimensional hypersurface for  $M$  3 dimensional 
  or 4 dimensional isotropic submanifold of $G_2$ manifold $M$. 
  
  We pullback the bundle $I$ on $T\Sigma^0$ via mapping 
   \begin{equation}
   Di: T^0(\Sigma) \rightarrow T^0(M)
   \end{equation}
 We end up with $Di^*(I)$. The claim is that $\Sigma$ defines a 
 Lagrangian subbundle in $Di^*(I)$. We pullback the tangent bundle $T\Sigma$ on $T^0(\Sigma)$ 
 (and denote it by the same letters). To avoid complicated notations 
 we denote the same letters the pullback of the tangent bundle $TM$ on 
 the compliment of the zero section in the tangent bundle $T^0(M)$.
 To define a Lagrangian subbundle we arrange
  the following mappings. 
 
 \begin{equation}
 s:T\Sigma \rightarrow Di^*(TM) \rightarrow Di^*(I) 
 \end{equation}  
 The Lagrangian subbundle $\cal{L}\Sigma $
  is the image of the composite map $s$.  
\section{Calculus of bubbles}
\subsection{The formulation of the problem}
Let us suppose  that we are in ${\Bbb R}^3$ we are given
 2 spheres in ${\Bbb R}^3$ with oriented normal bundles.  They intersect each other by a curve
or a union of curve from the space $\tilde Y$. The intersection carries a canonical orientation.  The
question is when it is possible to separate those two spheres.  In
the process of separating   them we want to keep spheres
intersecting by the curve in $\tilde Y$.  

\begin{Con}
 One can not separate two spheres,
 intersecting by a circle.  
\end{Con}
\begin{Thm} \label{T:ert}
Suppose that the process of separation is performed along good families of embeddings of spheres (see below), then one can not separate two spheres,
 intersecting by a circle. 
\end{Thm}

It is natural to formulate this problem in terms of
infinite dimensional  Lagrangian
intersections.  According to the previous sections
every embedded oriented 2-sphere defines a  Lagrangian
submanifold consisting of two components. We pick the one where curves oriented counterclockwise. It is clear that curves in the
intersections of spheres are exactly points of
intersections of corresponding Lagrangian manifolds.

Our aim is to construct an invariant, similar to the index of 
intersection which would be invariant when we continuously
 deform one of the Lagrangian submanifolds.

\subsection{Reminder of the basics of theory of finite
dimensional Lagrangian intersections}\label{S:qwr}

According to \cite{F} the story goes as follows.
In finite dimensions we are given
 a symplectic manifold $(M,\omega)$ and two connected 
 Lagrangian submanifolds $L_1,L_2$. There are two 
 definitions of intersection index. One is purely 
 homological, it is  the ordinary 
  index of intersection. The other is defined using 
  symplectic geometry. The former we want to 
  generalize in infinite dimensions.
  
Suppose that  Lagrangian submanifolds meet
 transversally at points $x,\dots,y$. 
 We want to define Maslov index $m(x,y)$.
  
  Define a space 
  
  \begin{equation}
  \begin{split}
  \Omega  =\Omega(L_,L_2) & =\\
  & =\{z \in Map([0,1],M)| z(0) \in L_1 \  \text{and} \  z(1) \in L_2 \}
  \end{split}
  \end{equation}
  
Constant maps in $\Omega$ correspond 
to points of intersection of $L_1$ and $L_2$. 
Define a path in $\Omega$ between $z_1$ 
and $z_2$ as a smooth map

\begin{equation}
u:[0,1]^2 \rightarrow M
\end{equation} 

so that $u(t,o)=z_1(t)$,$u(t,1)=z_2(t)$ and
 $u(0,t)\in L_1$ and $u(1,t)\in L_2$. 
 In particular if $z_1,z_2$ are constant,  $u$
  maps the boundary of $[0,1]^2$ 
  into $L_1 \cup L_2$ defining arcs in $L_1$ and 
  $L_2$ connecting $x=Im z_1$ and $y=Im z_2$. 
  In this situation Viterbo \cite{V} 
  defines an index by means of the Maslov class, i.e.,
   the generator
   
\begin{equation}
\mu \in H^1(\Lambda_n,{\Bbb Z})={\Bbb Z}
\end{equation}   
where $\Lambda_n \subset G({\Bbb R}^n,{\Bbb R}^{2n})$ 
is the set of Lagrangian n-planes in ${\Bbb R}^{2n}$.

 Consider a pullback of the tangent bundle $u^*TM$. 
 Since the base space $[0,1]^2$ is contractible, 
 we can define a trivialization
 \begin{equation}
 \Phi:u^*TM \rightarrow [0,1]^2 \times {\Bbb C}^n
 \end{equation}
 
  mapping the symplectic form on each fiber
 into standard form. Moreover, we can choose this trivialization
 so that it is constant on $\{0\}\times[0,1]$ and 
 on  $\{0\}\times[0,1]$, and so that 
 $\Phi(T_xL_1)=iT_xL_2$ and $\Phi(T_yL_1)=iT_yL_2$.
 Then define $m_u(x,y)$ by evaluating the 
 Maslov class over the closed loop
 \begin{equation}
 \partial [0,1]^2 \rightarrow \Lambda_n
 \end{equation}
 
 \begin{align}
 &(\tau,0) \rightarrow \Phi(T_{u(\tau,0)}L_1)\\
 &(1,t) \rightarrow e^{ {i\pi \over 2}  t}\Phi(T_{y}L_1)\\
 &(\tau,1) \rightarrow \Phi(T_{u(\tau,1)}L_2)\\
 &(0,t) \rightarrow e^{-{i \pi \over 2} t}\Phi(T_{x}L_2)
 \end{align}

\section{Relative(Maslov) index}

We want to define Maslov index which would be relevant to infinite dimensional setting.

 We are  given 2
embedded two (four) dimensional isotropic oriented submanifolds $\Sigma_1$
$\Sigma_2$, which are intersected transversally by two
or more curves lying in $\tilde Y$.  We pick two curves
$\gamma_1$ , $\gamma_2$ from the intersection, which
 belong to the same connected components of 
 $L(\Sigma_1)$ and $L(\Sigma_1)$.

 As in \ref{S:qwr} we pick two arcs , connecting
$\gamma_1,\gamma_2$ in $L(\Sigma_1)$ and
$L(\Sigma_2)$ respectively. We assume that there 
are no topological obstructions and there exist a map
\begin{equation}\label{E:vft}
u:[0,1]^2 \rightarrow \tilde Y
\end{equation} 
whose restriction $u(t,0)$ and $u(t,1)$ are the mentioned above arcs.
This map encodes a map of a handlebody

\begin{equation}
u^{\sharp}:[0,1]^2\times S^1 \rightarrow M
\end{equation}

$u(t,0,z)$ and $u(t,1,z)$ are maps of cylinders into
 $\Sigma_1$ and $\Sigma_2$ respectively, while $u(0,t,z)$
 and $u(1,t,z)$ do not depend on $t$.  As it was mentioned
 in \ref{S:asw} for every value $0 \leq a,b \leq 1$ we have
 a map

 \begin{equation}
u^{\sharp}(a,b): S^1 \rightarrow T^0(M), \ z \rightarrow u^{\sharp}(a,b,z)
\end{equation} 
and therefor
\begin{equation} \label{E:ert}
u^{\sharp} :[0,1]^2\times S^1 \rightarrow T^0(M)
\end{equation}

We pullback the bundle $I$ via $u^{\sharp}$ on
$[0,1]^2\times S^1$.  According to \ref {S:asw} in 3 and 7
dimensions bundle $I$ carries a complex structure.  This
being said it becomes clear that $(u^{\sharp})^*I$ is
trivial.  As in \ref{S:qwr} we choose it trivialization

\begin{equation}\label{E:xcv}
 \Phi:(u^{\sharp})^*I \rightarrow [0,1]^2\times S^1 \times {\Bbb C}^n \ n=1\ or \ 3
\end{equation}

mapping symplectic form on the fiber into standard form 
on ${\Bbb C}^n$. Moreover we chose this trivialization constant on 
$\{0\}\times[0,1]\times S^1$ and $\{1\}\times[0,1]\times S^1$.
 As we know (see \ref{S:tgu}) over $T^0(\Sigma)$ the bundle $I$
  has a Lagrangian subbundle $\cal{L}\Sigma$. Since the surfaces
  $\Sigma_1$ and $\Sigma_1$ are oriented then the restriction of each
   $\cal{L}\Sigma_1$ on $T^0(\Sigma)_1\cap T^0(\Sigma)_2$ is 
   trivial. We choose such trivialization $\Phi$ that 
\begin{equation} \label{E:fgh}  
   \Phi(\cal{L}(\Sigma_1))=i\Phi(\cal{L}(\Sigma_2))
\end{equation}
at  $T^0(\Sigma)_1\cap T^0(\Sigma)_2$.

Then we use the same formula as in \ref{S:qwr} to define a map:

\begin{equation}
 (\partial [0,1]^2)\times S^1 \rightarrow \Lambda_n \ n=1,3
\end{equation}

Observe that over points $[0,1]\times \{0\} \times S^1$ and
 $[0,1]\times \{1\} \times S^1$ of the base of the trivial
 bundle \ref{E:xcv} we have two well
 defined Lagrangian subbundles.  The mapping reads as :

\begin{align}
 &(\tau,0,z) \rightarrow \Phi(\cal{L}(\Sigma_1)_{u(\tau,0,z)})\\
 &(\tau,1,z) \rightarrow \Phi(\cal{L}(\Sigma_2)_{u(\tau,1,z)})
 \end{align}
 
 Over the points $\{0\}\times [0,1] \times S^1$ and 
 $\{1\}\times [0,1] \times S^1$ we want to employ the condition \ref{E:fgh}.
\begin{align}
&(1,t,z) \rightarrow e^{{i \pi  \over 2}t}\Phi(\cal{L}(\Sigma_1)_{\gamma_1})\\
 &(0,t,z) \rightarrow e^{-{i \pi \over 2}t}\Phi(\cal{L}(\Sigma_2)_{\gamma_2})
\end{align}

{\bf Observation} The pairing of the Maslov class $\mu$ with a 
distinguished 1 cycle 
\begin{equation}
\begin{split}
S=\{0\}\times\{0\}\times S^1\\
&\subset (\partial [0,1]^2)\times S^1
\end{split}
\end{equation}

 is is not necessarily zero.  Recall that we had a some
freedom in the choice of trivialization \ref{E:xcv}. Namely 
the trivializations are homotopically classified by 
$\pi_0(Map([0,1]^2\times S^1,\bold U(n))$, $n=1,3$. 
The last set is isomorphic to ${\Bbb Z}$. Since $\Sigma
_i$ ($i=1,2$) are oriented, we can always change 
trivialization  $\mu(S)$ to become zero.

Now we are able to give the definition. Pick some
cycle $P$, whose Poincar\'{e} pairing with $S$ is one.

\begin{Def}
The relative Maslov index of two curves 
$\gamma_i\subset \Sigma_1 \cap \Sigma_2$ ($i=1,2$) 
is defined as $m(\gamma_1,\gamma_2)=\mu(P)$
\end{Def} 

It should be clear that $m(\gamma_1,\gamma_2)$ doesn't 
depend on the  particular choice of $P$.

\subsection{Some properties of relative index}

We want to explore the dependence of  $m(\gamma_1,\gamma_2)$ 
on the choice of the map $u^{\sharp}$.

Now the dimension of the underlying manifold becomes essential.
\begin{Prop} \label{P:gyu}
If a manifold $M$ is 3-dimensional the index 
$m(\gamma_1,\gamma_2)$ is well defined $mod(4)$.
\end{Prop}

\begin{pf}
Suppose we have to maps $u_1$ and $u_2$ as in \ref{E:vft},
that match on the boundary $\partial [0,1]^2$. Together they define a map 
\begin{equation}
u_1\# u_2: S^2 \rightarrow \tilde Y
\end{equation}
Consider a pullback of a bundle 
$\tilde X \rightarrow \tilde Y$ on $S^2$ via the map $u_1\# u_2$. Denote the total space of the pullback by $W_{u_1\# u_2}$. Since $Diff^+(S^1)$ are homotopically equivalent to $S^1$, the space $W_{u_1\# u_2}$ is homotopically equivalent to either a lens space $S^3/ {\Bbb Z}_n$ or $S^2 \times S^1$. The equivalence is carried out by replacing a curve by a curve with natural parametrization. Again we encode the map $u_1\# u_2$, by a map:
\begin{equation} \label{E:yhu}
(u_1\# u_2)^{\sharp}:A \rightarrow M^3
\end{equation}
 $A$ stands for one of the lens space or $S^2 \times S^1$.  As in \ref{E:ert} we define a map , denoted by the same letter  as in \ref{E:yhu}, $(u_1\# u_2)^{\sharp}:A \rightarrow T^0(M^3)$. It is clear that the domain of this map sewed from the domains of $u_1^{\sharp}$ $u_2^{\sharp}$ along a torus.

Let  $I_{u_1\# u_2}$ be a pullback of $I$ on $A$. Observe $c_1(I_{u_1\# u_2})$ is given as an element $ \pi_0 ( \partial [0,1]^2 \times S^1, {\Bbb C}^*)$, which defines a difference between trivializations $\Phi_{u_1}$ and $\Phi_{u_2}$ over  $\partial [0,1]^2 \times S^1$. 

If $A$ is a lens space then $c_1(I_{u_1\# u_2})=0$ and there is no ambiguity in the definition of index. 

If $A$ is $S^2 \times S^1$ then $c_1(I_{u_1\# u_2})$ might be nonzero. According to  \ref{L:fvx} $c_1(I_{u_1\# u_2})$ is divisible by 2, another factor 2 due to the fact that natural map $U(n) \rightarrow \Lambda_n$ induces a multiplication by two in $\pi_1$.
\end{pf}
\begin{Prop}
If  $M$ is $G_2$ manifold the index 
$m(\gamma_1,\gamma_2)$ is well defined integer,
 depending on no choices.
\end{Prop}

\begin{pf}
The discussion goes as in the proof of \ref{P:gyu}. In a due moment we use \ref{L:wsd}. 
\end{pf}

\subsection{Definition of the invariant}

Let us suppose that we have two oriented surfaces $\Sigma_1$ $\Sigma_2$ which
intersect transversally by a link.  We chose a maximal set of
curves that lie in the same connected components of $L(\Sigma_1)$ and
$L(\Sigma_2)$.  We divide this set into two subsets:  two curves belong to the
same subset if their relative index is even and to different subset if it is
odd.  The difference of the cardinalities of these set is the invariant
$T(\Sigma_1,\Sigma_2,a)$.  By $a$ we denote the class of curves from the same
connected component of $L(\Sigma_1)$ and $L(\Sigma_2)$.  The invariant is
defined up to a sign.

We want to establish some kind of deformation invariance. It easy to construct 
a family of embedding of $\Sigma_1$ which violates conservation of 
$T(\Sigma_1,\Sigma_2,a)$. It is due to the fact that the total space
 $\tilde Y$ is not compact.

We discuss what are good deformations of a configuration $\Sigma_1$ $\Sigma_2$.
We keep an embedding $\phi : \Sigma_2 \rightarrow M$ fixed and vary the embedding of $\Sigma_1$ , $ \varphi: [0,1] \times \Sigma_1 \rightarrow M$ is the  family of embeddings. 

\begin{Def}
A family of embedding  $\varphi$ is considered as good
 if for every value of $t_0 \in [0,1]$ , $Im(\phi \cap Im(\varphi(t_0,.)$ is a link and $\phi$ and $\varphi$ are transversal.
\end{Def}

This definition implies hat there is a finite set of critical values of $t \in [0,1]$ for which $\phi$ and $\varphi(t,.)$ are not transversal. For these values the images still must intersect by a link.  As  a result we have a  simple local picture of critical intersection of two surfaces $\Sigma_1$ and $\Sigma_2$ along a knot $\gamma$. This configuration is diffeomorphic along $\gamma$ to the following one: rotate a parabola $z=x(x-2)+1$ about $z$ axis. The result is a surface $\Sigma_1'$ which is tangent to $x,y$ plane  (surface $\Sigma_2'$) at the points of intersection, which is a circle $x^2+y^2=1$ (a knot $\gamma'$).
Then we have a simple 
\begin{Prop}
The invariant $T(\Sigma_1,\Sigma_2,a)$ is preserved along good deformations.
\end{Prop}
\begin{pf}
It is clear that we can suppose that interval $[0,1]$ contains on;y one critical value of  of parameter $t$ for $\varphi$. The the direct computation shows that a pair of knots which annihilates at this critical value has odd relative index.
\end{pf}
 The proof of the theorem \ref{T:ert} now becomes obvious.
It seems  reasonable that to prove a similar theorem when two spheres intersect each other by large number of circles one has to develop some version of Floer complex. 

\section{Complex structure } 
We would like to go a little bit further and exhibit a
quasi-complex structure, which exists on $\tilde Y$.

In 3 dimensions it was made in \cite{Bry}. We want to do a 
similar thing for $G_2$ manifolds. According to section \ref{S:ert}
 we can  canonically associate to the 3-form $\lambda$ a 
 Riemannian metric $g$ and $(2,1)$ tensor $[.]$. The tensor
  $[.]$ is an analog of the vector product in 3 dimensions. This being 
  said the definition  of the complex structure becomes a rephrasing of the one in  \cite{Bry}:we identify the tangent space to a
   point of $\tilde Y$, represented by a loop $\gamma$ with a 
   normal bundle to $\gamma$. Plugging tangent $\overset{.}{\gamma}$ into $[,]$ ,
   get an operator $[\overset{.}{\gamma},.]$, acting in the normal bundle.
    The equation \ref{E:gty}  secure that it square is $-id$.
It is unclear to me if $d^*$ closeness of $\lambda$ implies integrability of the complex structure. It is more likely that it is always nonintegrable. Here are some evidence for it. Consider a simplest example $M^7={\Bbb R}^7$. Instead of  unparametrized space of loops take the  unparametrized space of paths originating at the origin. This space carries a quasi-complex structure similar to the one described.  There is a complex map of $S^6$(with a classical quasi-complex structure) to the space of paths.  A point on the sphere generate a ray through it. So the complex structure on the space of paths must be nonintegrable.  
  
\section{One Hamiltonian system}

On the space $\tilde Y$ associated with ${\Bbb R}^7$ there is a canonical function-the length of the curve. We want to understand better the structure of the corresponding Hamiltonian flow. It is known that the analog of this system in 3 dimensions is integrable. We convinced that in 7 dimensions integrability still holds, though we where unable to prove it. Here are some partial results
\begin{Prop}
The system ha s three first integrals.
\end{Prop}
\begin{pf}
Following \cite{Bry} the Hamilton equation have form

\begin{equation}\label{E:iop}
\overset{.}{\gamma(z,t)}=[\gamma'(z,t),\gamma(z,t)'']
\end{equation}

We use the notations; $z$ is the parameter on the circle, $t$ stands for time, $.$ partial derivative with respect to $t$ , $'$ partial derivative with respect to $x$.Introduce a new function $T=\gamma'$. Differentiating \ref{E:iop} with respect to $z$ we get 

\begin{equation}
\overset{.}{T}=[T,T'']
\end{equation}

Apriori the value of the integral with a density $(\gamma',\gamma')$ is conserved,

Then we have an equality :
\begin{equation}
1/2(T',T')^.=-\lambda(T,T',T'')''
\end{equation}
So $1/2(T',T')$ leads to a conservation law.
The last density we have in mind is $\lambda(T,T',T'')$, it also leads to a conservation law. Of course we anticipate the recursion of the kind $\overset{.}{P_i}=P_{i+1}'$ ($P_i$ are conservation law) , but we were unable to check it in full generality. Unfortunately the method which is used to prove a similar fact in 3 dimensions uses heavily Lie algebra structure of the vector product. 
 \end{pf}

\end{document}